\documentclass[12pt]{amsart}
\usepackage{times}
\usepackage{amsmath}
\usepackage{amsfonts}
\usepackage{amssymb}
\newcommand{\dis}{\displaystyle}

\newcommand{\Rr}{\mathbb{R}}

\newcommand{\Qq}{\mathbb{Q}}

\newcommand{\Zz}{\mathbb{Z}}

\newtheorem{theorem}{Theorem}
\newtheorem{lemma}[theorem]{Lemma}

\newcommand{\beq}{\begin{eqnarray}}
\newcommand{\eeq}{\end{eqnarray}}
\newcommand{\beqn}{\begin{eqnarray*}}
\newcommand{\eeqn}{\end{eqnarray*}}

\begin{document}
\title[Heegner Divisors in the Period Space of Enriques Surfaces]{Irreducible Heegner Divisors in the Period Space of Enriques Surfaces}
\author{Caner Koca}
\address{Stony Brook University \\ Department of Mathematics \\
Stony Brook, NY 11794-3651, USA}
\email{caner@math.sunysb.edu \\ {\it Fax:} 1-631-632 7631}
\author{Al\.i S\.inan Sert\"{o}z}
\address{B\.ilkent University \\ Department of Mathematics \\
TR-06800 Ankara, Turkey}
\email{sertoz@bilkent.edu.tr \\ {\it Fax:} 90-312-266 4579}
\renewcommand{\subjclassname}{\textup{2000} Mathematics Subject Classification}
\begin{abstract}
Reflection walls of certain primitive vectors in the anti-invariant sublattice of the K3 lattice define Heegner divisors in the period space of Enriques surfaces. We show that depending on the norm of these primitive vectors, these Heegner divisors are either irreducible or have two irreducible components. The two components are obtained as the walls orthogonal to primitive vectors of the same norm but of different type as ordinary or characteristic.
\end{abstract}
\maketitle


\section{Introduction}
The deck transformation of the universal covering of an Enriques surface
induces an involution on the second integral cohomology of the cover, which in turn,
under a marking, defines an involution $\tau$ on the K3 lattice $\Lambda$. If $\Lambda^-$ denotes the eigenspace of $\tau$ corresponding to the eigenvalue $-1$, we define the period domain of Enriques surfaces as $\mathcal{D}=\{ [v]\in \mathbb{P}(\Lambda^-\otimes \mathbb{C})| \langle v,v\rangle=0, \; \langle v,\overline{v}\rangle >0\}$, where $\langle\cdot ,\cdot \rangle$ denotes the bilinear form of $\Lambda$. If $O(\Lambda^-)$ denotes the orthogonal group of $\Lambda^-$, then $\mathcal{D}/O(\Lambda^-)$ is the period space of Enriques surfaces.

We define $H_{\ell}=\{ [v]\in \mathbb{P}(\Lambda^-\otimes \mathbb{C})| <v,\ell>=0 \}$, for all primitive $\ell\in\Lambda^-$. For any integer $n$, let $\mathcal{D}_n$ denote $\bigcup H_{\ell}$ where the union is taken over all primitive $\ell\in \Lambda^-$ with $\langle \ell,\ell\rangle=2n$. The Heegner divisor $\mathcal{H}_n$  is now defined as $\mathcal{D}_n/O(\Lambda^-)$.

We show that $\mathcal{H}_n$ is irreducible for odd $n$, and has two
disjoint irreducible components for even $n$. Moreover we have a full description of
these components when $n$ is even. For this recall that an element $\ell\in \Lambda^-$ is
called characteristic if $\langle \ell,\nu\rangle \equiv \langle \nu,\nu\rangle \; {\rm
mod} \; 2$ for all $\nu\in \Lambda^-$, and is called ordinary otherwise. Let $\mathcal{D}_n^c$
(resp. $\mathcal{D}_n^o$) denote $\bigcup H_{\ell}$ where
the union is taken over all primitive characteristic (resp. ordinary) $\ell\in \Lambda^-$ with
$\langle \ell,\ell\rangle=2n$. Define the two Heegner divisors $\mathcal{H}_n^c$ and
$\mathcal{H}_n^o$ as $\mathcal{D}_n^c/O(\Lambda^-)$ and $\mathcal{D}_n^o/O(\Lambda^-)$ respectively. It follows that $\mathcal{H}_n=\mathcal{H}_n^c \bigsqcup\mathcal{H}_n^o$ for even $n$.

For this we first count the number of orbits in $\Lambda^-$ under
the action of its isometries, see Theorem~\ref{theorem:1}, and use this to describe Heegner divisors in section~\ref{last}. As a byproduct of our proof we obtain a simple description of even type vectors of $\Lambda^-$ in section~\ref{penultimate}.


\section{Definitions and the main theorem}

We consider an Enriques surface $S$ whose K3 cover we denote by
$\tilde{S}$. There is a fixed point free involution $\iota :
\tilde{S}\longrightarrow\tilde{S}$ which extends to an involution
on cohomology $\tau:H^2(\tilde{S},\Zz)\longrightarrow
H^2(\tilde{S},\Zz)$. Combined with a marking $\phi
:H^2(\tilde{S},\Zz)\longrightarrow \Lambda$ we can consider
$\tau$ as acting on $\Lambda$, where $\Lambda$ is the K3
lattice $E_8^2\oplus U^3$. Here $E_8$ is the unimodular, negative
definite even lattice of rank 8, and $U$ is the hyperbolic plane.
We denote the $(-1)$-eigenspace of $\tau$ in $\Lambda$ by
$\Lambda^-$. It is known that $\Lambda^-\cong E_8(2)\oplus
U(2)\oplus U$. We will work with a fixed marking and take
$\Lambda^-=E_8(2)\oplus U(2)\oplus U$. We denote the orthogonal
group of $\Lambda^-$, the group of its self isometries, by
$O(\Lambda^-)$.

The lattice product of two elements $\omega_1$, $\omega_2$ in any
lattice $L$ will be denoted by $\langle \omega_1,\omega_2\rangle _L$ where we
can drop the subscript $L$ if it is clear from the context. We say
that $\omega$ is an $m$-vector, or an element of norm $m$ in $L$
if $\langle \omega,\omega\rangle _L=m$. A nonzero element $\omega$ in a lattice
$L$ is called primitive if $\frac{1}{n}\omega$ is not an element
in $L$ for any nonzero integer $n\neq \pm 1$.

For any integer $n$ we use $(n)$ to denote the lattice $\Zz\cdot
e$ where $\langle e,e\rangle =n$. For any lattice $L$ and any positive integer
$n$, $L^n$ denotes the direct sum of $n$ copies of $L$. For any
nonzero integer $n$, $L(n)$ denotes the same $\Zz$-module $L$
except that the inner product is modified by $n$ as
$\langle \omega_1,\omega_2\rangle _{L(n)}=n\langle \omega_1,\omega_2\rangle _{L}$ for any
$\omega_1,\; \omega_2\in L$. We denote by $I_{s,t}$ the odd
unimodular lattice $(1)^s\oplus (-1)^t$, where $s$ and $t$ are
nonnegative integers.

We define the type of an element $\omega\in L$, for any lattice
$L$, to be  characteristic  if
$\langle \omega,\eta\rangle \equiv \langle \eta,\eta\rangle  \mod 2 \; \; \mbox{for all~} \eta\in L,$
and ordinary otherwise.   Furthermore when $L$ is unimodular and $\omega\in L$ is characteristic, then $\langle \omega,\omega\rangle \equiv n-m \mod 8$, where $(n,m)$ is the signature of $L$, see
\cite[Theorem 4]{wall}, \cite[II.5.2]{milnor2}.

\begin{theorem} \label{theorem:1}
If $B$ is an even, unimodular, indefinite lattice, then the number of distinct orbits of primitive $(2n)$-vectors
in $B(2)\oplus U$ under the action of $O(B(2)\oplus U)$ is $1$ when $n$ is odd, and $2$ when $n$ is even.
\end{theorem}

Note that since $\Lambda^-=(E_8\oplus U)(2)\oplus U$, this theorem applies to
$\Lambda^-$. The proof of the theorem will be given in section~\ref{proof}, and its implementation to Heegner divisors in the period space of Enriques surfaces will be given in section~\ref{last}.

We remind that some special cases of the theorem for $\Lambda^-$
were proven by Namikawa (for $n=-1,-2$, see \cite{namikawa}
Theorems 2.13, 2.15),  by Sterk (for $n=0, -2$, see \cite{sterk}), and by Allcock (for $n=0,-1,-2$, see
\cite{allcock}) .


\section{A crucial dilatation}

In order to count the orbits of $O(\Lambda^-)$ in $\Lambda^-$ we will use a dilatation which can be deduced from the following setup.  If $A=B(2)\oplus U$, then $\dis \left(
(2^{-1/2})A\right)^\ast$ is isomorphic to $B\oplus U(2)$ and this
lies isometrically in $B\oplus I_{1,1}$. If the signature of $B$
is $(s-1,t-1)$, then $B\oplus I_{1,1}\cong  I_{s,t}$. It can be observed that
$O(A)\cong O(I_{s,t})$, see
\cite{allcock,Kondo}.

We will construct a dilatation of $B(2)\oplus U$ into
$(B\oplus I_{1,1})\otimes \Qq$ through which we  recover the isomorphism
of the above orthogonal groups. This has the added advantage of giving explicit description of the components of the Heegner divisors under consideration. We start with a lemma.

\begin{lemma} \label{lemma2}
An element $\omega=(a_1,\dots,a_r,m,n)\in B\oplus I_{1,1}$ is
characteristic if and only if $m$ and $n$ are odd and all $a_i$'s
are even.
\end{lemma}

\begin{proof}
One way is straightforward. So take $\omega$ to be characteristic.
Take $\eta=(0,\dots,0,1,0)$. Since $\langle \omega,\eta\rangle =m$ and
$\langle \eta,\eta\rangle =1$, it follows that $m$ is odd. Similarly we can show $n$ to be odd also.
Now take $\eta=(c_1^i,\dots,c_r^i,0,0)$ where
$\vec{c}^{~i}$ is the $i$-th row of $B^{-1}$, $i=1,\dots,r$.
If $\vec{a}$ denotes the column vector with $a_i$ in the
$i$-th place, $i=1,\dots,r$, then
we have
\[ \langle \eta, \omega\rangle=(c_1^i,\dots,c_r^i,0,0)\left( \begin{array}{ccc}
B & 0 & 0 \\
0 & 1 & 0 \\
0 & 0 & -1
\end{array} \right) \left( \begin{array}{c} a_1 \\ \vdots \\ a_r \\ m \\ n \end{array} \right) = \vec{c}^{~i} B \vec{a} = a_i, \; \; i=1,\dots,r. \]
On the other hand we have $\langle \eta,\eta\rangle \equiv 0 \mod 2$ since
$B$ is even. This forces all $a_i$'s to be even since $\omega$ is assumed to be characteristic.
\end{proof}

For any lattice $L$ of rank $n$, fix a basis and  define the following
$\Zz$-module:
\begin{eqnarray*}
( \frac{1}{2}) L &=& \{ \frac{1}{2}\omega\in \Rr^n\; | \;
\omega\in L\; \mbox{and coordinates of } \omega \mbox{~are} \\ &&
\quad \quad \quad \mbox{~either all odd or all even}. \; \}.
\end{eqnarray*}
Extend the inner product of $L$ to this $\Zz$-module linearly to give it a lattice structure. This description is dependent on the chosen basis but it  is a particular intermediate construction we use below to obtain  coordinate free isomorphisms.

\begin{lemma} $O(B\oplus ( \frac{1}{2})I_{1,1}) \cong O(B\oplus I_{1,1})$.
\end{lemma}
\begin{proof}
$B\oplus I_{1,1}$ is a submodule of $B\oplus (
\frac{1}{2})I_{1,1}$, where both are considered as submodules of
$\Zz^{r+2}$. It can be shown that with respect to some basis of
$B\oplus ( \frac{1}{2})I_{1,1}$ which uses the standard basis in
$I_{1,1}$, any isometry of $B\oplus ( \frac{1}{2})I_{1,1}$ is
represented by an integral matrix which then defines an isometry
on $B\oplus I_{1,1}$. This defines an injective morphism from
$O(B\oplus ( \frac{1}{2})I_{1,1})$ to $O(B\oplus I_{1,1})$. Now
take any isometry $g$ in $O(B\oplus I_{1,1})$ and
$\omega=(a_1,\dots,a_r,m+\frac{1}{2},n+\frac{1}{2})\in B\oplus
(\frac{1}{2})I_{1,1}$. It follows that $2\omega$ is a characteristic
element in $B\oplus I_{1,1}$ by lemma \ref{lemma2}, and $g(2\omega)$ is
characteristic, see \cite{wall}. Finally by the same lemma
$\frac{1}{2}g(2\omega)$ is an element of $B\oplus
(\frac{1}{2})I_{1,1}$ and thus $g$ defines an element in $O(B\oplus ( \frac{1}{2})I_{1,1})$.
\end{proof}

In the following lemma we construct the crucial dilatation $\phi$.

\begin{lemma}\label{iso}
If $B$ is an even, unimodular lattice of signature $(s-1,t-1)$,
then $O(B(2)\oplus U)\cong O(B\oplus I_{1,1})\cong O(I_{s,t})$.
\end{lemma}

\begin{proof}
Let $e_1,\dots,e_r$ be a basis of $B(2)$, where $r=s+t-2$. Let
$u,v$ be a basis of $U$ such that $\langle u,u\rangle =\langle v,v\rangle =0$ and $\langle u,v\rangle =1$.
Any element of $B(2)\oplus U$ is of the form
$a_1e_1+\cdots+a_re_r+b_1u+b_2v$ for some integers $a_i$ and
$b_i$. We denote this element, as above, with the vector
$(a_1,\dots,a_r,b_1,b_2)$. Similarly for the lattice $B\oplus
(\frac{1}{2})I_{1,1}$ where now $x,y$ are basis for $I_{1,1}$ with
$\langle x,x\rangle =-\langle y,y\rangle =1$ and $\langle x,y\rangle =0$. We define a map
\begin{eqnarray*}
\phi:B(2)\oplus U & \rightarrow & B\oplus (\frac{1}{2})I_{1,1} \\
(a_1,\dots,a_r,b_1,b_2) & \mapsto &
(a_1,\dots,a_r,\frac{b_1+b_2}{2},\frac{b_1-b_2}{2})
\end{eqnarray*}

This is a $\Zz$-module isomorphism with
\[ \langle \omega_1,\omega_2\rangle =2\langle \phi(\omega_1),\phi(\omega_2)\rangle , \; \; \mbox{for~~}
\omega_i\in B(2)\oplus U. \] At this point it is clear that
$O(B(2)\oplus U)\cong O(B\oplus (\frac{1}{2})I_{1,1})$. The
previous lemma gives the isomorphism to $O(B\oplus I_{1,1})$, which is in turn isomorphic to  $O(I_{s,t})$ via the classification theory of unimodular odd
lattices.  \end{proof}


\section{Counting the orbits}\label{proof}

In this section we give a constructive proof of theorem \ref{theorem:1}.
In what follows $B$ is an even, unimodular, indefinite matrix of signature
$(s-1,t-1)$.

Lemma \ref{iso} reduces the problem of counting
orbits in $B(2)\oplus U$ to the same problem in $I_{s,t}$. The
number of orbits in $I_{s,t}$ is known, \cite[Theorem 4]{wall}:

\begin{theorem}[C.T.C. Wall]
If $s,t$ are each at least 2, then $O(I_{s,t})$ acts transitively
on primitive vectors of given norm and type (i.e. characteristic
or ordinary). Moreover, if a vector is characteristic then its
norm is congruent to $s-t \mod 8$. \hfill $\Box$
\end{theorem}

Since $B$ is even, indefinite and unimodular, we necessarily have
$s-t\equiv 0 \mod 8$, see \cite{milnor}.

\begin{proof}[Proof of theorem \ref{theorem:1}] \mbox{} \\
{\bf Case 1: $n$ is an odd integer. }

Let $\omega=(a_1,\dots,a_r,b_1,b_2)$ be  a primitive vector in
$B(2)\oplus U$.  We have $\langle \omega,\omega\rangle =4k+2b_1b_2=2n$ with $n$
odd, so $b_1$ and $b_2$ are odd. Then
$\phi(\omega)=(a_1,\dots,a_r,\frac{b_1+b_2}{2},\frac{b_1-b_2}{2})$
is integral, primitive and since $n \not\equiv 0 \mod 8$ it is
ordinary. Since $O(I_{s,t}) \cong O(B\oplus I_{1,1})$, all such
elements are in the same orbit by Wall's theorem above. And since
$O(B(2)\oplus U)\cong O(B\oplus I_{1,1})$, all primitive elements
of norm $2n$ lie in the same orbit in $B(2)\oplus U$.
This orbit is not empty since
$\omega=(0,\dots,0,1,n)$ lies in it.

{\bf Case 2: $n$ is an even integer. }

$\langle \omega,\omega\rangle =4k+2b_1b_2=2n\equiv 0 \mod 4$. In this case $b_1$
and $b_2$ cannot both be odd.

{\bf Case 2.1: Only one of $b_1$  and $b_2$ is even.}

In this case
$\phi(\omega)=(a_1,\dots,a_r,\frac{b_1+b_2}{2},\frac{b_1-b_2}{2})$ and is fractional. But
$2\phi(\omega)$ is integral, primitive and by lemma \ref{lemma2}
it is characteristic. To show that this orbit is nonempty take
$\omega=(0,\dots,0,1,n)$ in $B(2)\oplus U$. It
is primitive, of norm $2n$ with only one of $b_1$ and $b_2$ even.
So this case contributes one orbit.

{\bf Case 2.2 $b_1$ and $b_2$ are both even.}

In this case $\phi(\omega)$ is  integral, primitive and ordinary.
Hence this case contributes one orbit if we can show the existence
of a primitive vector $\omega\in B(2)\oplus U$ of norm $2n$, with
both $b_1$ and $b_2$ even.

Since $B$ is even, unimodular and indefinite,  $B\cong U^i\oplus
E_8(\pm 1)^j$ where  $i > 0$. For any integer $k$, $(1,k)\in U(2)$
is a primitive element of norm $4k$, so in particular $B(2)$
contains a primitive element $\omega'$ of norm $4k$. Let $n=2k$.
Then $\omega=\omega'+(0,0)\in B(2)\oplus U$ is primitive of norm
$2n$ with $b_1$ and $b_2$ even. So this orbit is not empty either.

Both cases 2.1 and 2.2 contribute one orbit each, so when $n$ is
even  there are two orbits.
\end{proof}

We remark in passing that Sterk uses the isotropic vectors $e$,
$e'$, $e'+f'+\omega$, $e'+2f'+\alpha$ and $2e+2f+\alpha$, for
notation see \cite[propositions 4.2.3 and 4.5]{sterk}. He shows
that  $e$ belongs to an orbit disjoint from the orbit to which all
the others belong. In our point of view, $\phi(e)$ is fractional
with $2\phi(e)$ characteristic whereas the $\phi$ image of the
others are integral and ordinary, which explains the existence of
two distinct orbits.

\section{Even type vectors}\label{penultimate}

Regarding $\Lambda^+,\Lambda^-$ as sublattices of $\Lambda$, a
primitive $(2n)$-vector $\omega\in\Lambda^-$ is defined to be
\textit{even} if there is a primitive $(2n)$-vector
$\omega'\in\Lambda^+$ such that $\omega+\omega'\in 2\Lambda$,
(\textit{cf}. \cite[Proposition 2.16]{namikawa}). We give an explicit description of even type primitive vectors in terms of the dilatation $\phi$.

\begin{theorem}
Let $n$ be an even integer. A primitive $(2n)$-vector
$\omega\in\Lambda^-$ is of even type if and only if its image
$\phi(\omega)$ has integral coordinates.
\end{theorem}
\begin{proof}
Since any self-isometry of $\Lambda^\pm$ extends to a
self-isometry of $\Lambda$ (see \cite[Theorem 1.4]{namikawa})
without loss of generality we can fix a primitive embedding of
$E_8(2)\oplus U(2)\oplus U$ into $E_8^2\oplus U^3$ and prove the
statement for the image of this embedding. Thus, we fix the
following embedding
\begin{eqnarray*}
E_8(2)\oplus U(2)\oplus U & \hookrightarrow & E_8^2\oplus U^3\\
(e,u,v)   & \mapsto         & (e,-e,u,-u,v)
\end{eqnarray*}
and identify the domain with its image. The orthogonal complement
of the image is precisely the image of the primitive embedding
\begin{eqnarray*}
E_8(2)\oplus U(2)         & \hookrightarrow & E_8^2\oplus U^3\\
(e,u)     & \mapsto           & (e,e,u,u,0).
\end{eqnarray*}
Moreover, the primitive $(2n)$-vectors $\omega\in\Lambda^-$ with
integral images are transitive by Case 2.2. Therefore, it suffices
to prove the statement for particular vectors:

Let $\omega=(0,\dots,0,k,1,0,0)\in E_8(2)\oplus U(2)\oplus U$ be a
primitive vector with $\omega^2=2n=4k$, and identify $\omega$ with
its image in $E_8^2\oplus U^3$ with coordinates
$(0,\dots,0,k,1,-k,-1,0,0)$, by the above embedding. Notice that
$\phi(\omega)$ has integral coordinates. Now choose
$\omega'=(0,\dots,0,k,1)\in E_8(2)\oplus U(2)$ which corresponds
similarly to $(0,\dots,0,k,1,k,1,0,0)\in E_8^2\oplus U^3$. Now it
is clear that $\omega+\omega'\in 2 ( E_8^2\oplus U^3)$.

On the other hand, the $\phi$ image of the vector
$\omega=(0,\dots,0,2k,1)\in E_8(2)\oplus U(2)\oplus U$ has
fractional coordinates, and for no vector $\omega'$ in
$E_8(2)\oplus U(2)$ can we have $\omega+\omega'\in 2 (E_8^2\oplus
U^3)$, because the last coordinate of $\omega+\omega'$ is always
1.
\end{proof}

\section{Heegner Divisors}\label{last}

Here we discuss the irreducibility of certain Heegner divisors in the period space of Enriques surfaces. The period domain is defined as
\[ \mathcal{D}=\{ [v]\in \mathbb{P}(\Lambda^-\otimes \mathbb{C})|
\langle v,v\rangle=0, \; \langle v,\overline{v}\rangle > 0\} \]
and the period space is then $\displaystyle \mathcal{D}/O(\Lambda^-)$, see \cite{barth}.
For any primitive  $\ell\in\Lambda^-$ define
\[ H_{\ell}=\{ [v]\in \mathbb{P}(\Lambda^-\otimes \mathbb{C})| \langle v,\ell\rangle =0 \}, \]
which is the hyperplane orthogonal to $\ell$. For any integer $n$ define
\[ \mathcal{D}_n= \bigcup H_{\ell} \]
where the union is taken over all primitive $\ell\in \Lambda^-$ with $\langle \ell,\ell\rangle=2n$. The Heegner divisor $\mathcal{H}_n$ in the period space is now defined as
\[ \mathcal{H}_n =\mathcal{D}_n/O(\Lambda^-). \]
It follows from theorem~\ref{theorem:1} that $\mathcal{H}_n$ is irreducible for odd $n$, and has exactly two disjoint irreducible components for even $n$. Moreover the proof of theorem~\ref{theorem:1} gives a full description of the components of $\mathcal{H}_n$ for even $n$. For this define the two Heegner divisors
\[ \mathcal{H}_n^c \; \left( \mbox{resp.~}\mathcal{H}_n^o\right) =\left(  \bigcup H_{\ell}\right)/O(\Lambda^-), \]
where the union is taken over all primitive characteristic (resp. ordinary) $\ell\in \Lambda^-$ with
$\langle \ell,\ell\rangle=2n$.
Then we have
\begin{theorem}
\[ \mathcal{H}_n=\mathcal{H}_n^c \bigsqcup\mathcal{H}_n^o, \mbox{~for even~} n. \]
\end{theorem}
\mbox{}\hfill $\Box$

It is known that the period of an Enriques surface $Y$ lies in $H_{-2}^o$ if and only if $Y$ contains a nonsingular rational curve, see \cite[Theorem 6.4]{namikawa}. It will be interesting to find similar geometric descriptions for the other Heegner divisors.
Further characterizations of nodal Enriques surfaces were also obtained by Cossec \cite{cossec}.  For a recent investigation of this problem we refer to \cite{artebani}.


{\bf Acknowledgements: } We thank I. Dolgachev and D. Allcock for
numerous correspondences which helped us to clarify our arguments.
We also thank our colleagues A. Degtyarev, A. Klyachko and E.
Yal\c{c}{\i}n for several comments.



\newpage


\begin{thebibliography}{99}



\bibitem{allcock} Allcock, D., The period lattice for Enriques surfaces,\\
Math. Ann., 317 (2000), 483-488.

\bibitem{artebani} Artebani, M., Heegner divisors in the moduli space of genus three curves, {\tt arXiv.math.AG/0510199}.

\bibitem{barth} Barth, W. P., Hulek, K., Peters, C. A. M., Van De Ven, A., {\it Compact Complex Surfaces}, Second Enlarged Edition, Springer-Verlag, 2004.

\bibitem{cossec} Cossec, F., On the Picard group of Enriques surfaces, Math. Ann., 271 (1985), 577-600.

\bibitem{Kondo} Kond\={o}, S., The rationality of the moduli space of
Enriques surfaces, \\
Compositio Math., 91 (1994), 159-173.


\bibitem{milnor} Milnor, J., On simply connected 4-manifolds, \\
Symposium Internacional de Topologia Algebraica, pp 122-128, 1958.

\bibitem{milnor2} Milnor, J., Husemoller, D., Symmetric Bilinear Forms, Springer-Verlag, 1973.

\bibitem{namikawa} Namikawa, Y., Periods of Enriques surfaces, \\
Math. Ann., 270 (1985), 201-222.



\bibitem{nikulin} Nikulin, V.V., Integral Symmetric Bilinear Forms and Some
of their
Applications,\\
Math. USSR Izvestija, Vol. 14, No.1, (1980), 103-167.






\bibitem{sterk} Sterk, H., Compactifications of the period space of Enriques
surfaces, I,\\ Math. Z., 207 (1991), 1-36.


\bibitem{wall} Wall, C.T.C, On the orthogonal groups of unimodular quadratic
forms, \\
Math. Ann., 147 (1962), 328-338.





\end{thebibliography}
\end{document}